\documentclass[a4paper,10pt]{article}
\usepackage[utf8]{inputenc}
\usepackage{graphicx,amsmath,amsfonts,amssymb,graphbox}

\title{ Model Reduction Using Sparse Polynomial Interpolation for the Incompressible Navier-Stokes Equations }
\author{Martin W. Hess, Gianluigi Rozza}

\def\cl {\nonumber \\}
\def\el {\nonumber } 

\begin{document}

\maketitle

\begin{abstract}
This work investigates the use of sparse polynomial interpolation as a model order reduction method for the incompressible Navier-Stokes equations.
Numerical results are presented underscoring the validity of sparse polynomial approximations and comparing with established reduced basis techniques.
Two numerical models serve to access the accuracy of the reduced order models (ROMs), in particular parametric nonlinearities arising from curved geometries are investigated in detail.
Besides the accuracy of the ROMs, other important features of the method are covered, such as offline-online splitting, run time and ease of implementation.
The findings establish sparse polynomial interpolation as another instrument in the toolbox of methods for breaking the curse of dimensionality.
\end{abstract}

\section{Introduction}

Reduced order methods compute a low-order approximation to solutions of a parametrized partial differential equation (PDE) or derived quantities of interest of the solutions.
The reduced order model (ROM) evaluates the low-order approximation in a fast and computationally efficient manner, such that an inexpensive yet accurate approximate solution is available.
For an introduction and overview to ROM methods, in particular the reduced basis (RB) method, 
see for example \cite{hesthaven2015certified}, \cite{QuarteroniManzoniNegri}, \cite{doi:10.1137/1.9781611974829}, \cite{HBMOR_vol1}, \cite{HBMOR_vol2}, \cite{HBMOR_vol3}.

This work aims to establish sparse polynomial interpolation as a ROM method for the incompressible Navier-Stokes equations.
Two numerical models serve to access the accuracy of the low-order approximation computed from the sparse polynomial interpolation method.
More specifically, both models have previous results available obtained with the RB method \cite{10.1007/978-3-030-39647-3_45}, \cite{MR4099821}, which allows not only to compare the accuracy between both methods, but also
the run time, implementation effort and other desirable features, such as offline-online splitting.

There exists a significant body of literature on sparse polynomial interpolation, see \cite{10.1007/s10208-013-9154-z}, \cite{CHKIFA2015400}, \cite{doi:10.1142/S0219530511001728} and \cite{MR3780742}, which 
establishes the theoretical rationale behind the method and the relevant algorithms in an abstract setting. 
In particular, bounds on the accuracy of the ROM approximations w.r.t. the full order model (FOM) are established under holomorphy assumptions on the parameter-to-solutions mappings.
This could be a promising step towards breaking the so-called \emph{curse of dimensionality}.
The \emph{curse of dimensionality} refers to a (sub-)exponentially increasing computational effort with increasing parameter space dimension.
That potentially renders all problems unfeasible as the resolution in parameter space is increased in uncertain models for instance, but also poses a bottleneck for complex applications where many parameters are present.

However, this work aims to establish the sparse polynomial interpolation as a numerical method in incompressible computational fluid dynamics (CFD) and thus focuses more on the 
comparison to the already established RB method. Two non-trivial models are investigated. The first model has one parametric variation in geometry, which is affinely parametrized. 
The second model has two parametric variations, the kinematic viscosity and the curvature. The curvature introduces a nonlinear parameter dependency, which can not be easily resolved by RB methods.

The remainder of the work is structured as follows. 
Section 2 introduces the incompressible Navier-Stokes equations and the non-linear solver, while section 3 recapitulates the sparse interpolation procedure.
Section 4 provides and discusses the numerical results and section 5 concludes the findings and gives a brief outlook.

\section{Model Setup}

Let $\Omega \in \mathbb{R}^2$ be the computational domain.
Incompressible, viscous fluid motion in spatial domain $\Omega$ over a time interval $(0, T)$
is governed by the incompressible
 \emph{Navier-Stokes} equations:
%\eqref{Hess:NSE0} - \eqref{Hess:NSE1}:
\begin{eqnarray}
\frac{\partial \mathbf{u}}{\partial t} + \mathbf{u} \cdot \nabla \mathbf{u} &=& - \nabla p + \nu_{\text{visc}} \Delta \mathbf{u} + \mathbf{f}, \label{Hess:NSE0}  \\
\nabla \cdot \mathbf{u} &=& 0, 
\label{Hess:NSE1}
\end{eqnarray}
where $\mathbf{u}$ is the vector-valued velocity, $p$ is the scalar-valued
pressure, $\nu_{\text{visc}}$ is the kinematic viscosity and $\mathbf{f}$ is a body forcing.
\noindent Boundary and initial conditions are prescribed as 
\begin{eqnarray}
\mathbf{u} &=& \mathbf{d} \quad \text{ on } \Gamma_D \times (0, T), \\
\nabla \mathbf{u} \cdot \mathbf{n} &=& \mathbf{g} \quad \text{ on } \Gamma_N \times (0, T), \\
\mathbf{u} &=& \mathbf{u}_0 \quad \text{ in } \Omega \times 0,
\label{Hess:NSE_boundaryCond}
\end{eqnarray}

\noindent with $\mathbf{d}$, $\mathbf{g}$ and $\mathbf{u}_0$ given and $\partial \Omega = \Gamma_D \cup \Gamma_N$, $\Gamma_D \cap \Gamma_N = \emptyset$.
The \emph{Reynolds} number $Re$, which characterizes the flow (\cite{Holmes:2012}), depends on $\nu_{\text{visc}}$,
a characteristic velocity $U$, and a characteristic length $L$:
\begin{equation}\label{eq:re}
Re = \frac{UL}{\nu_{\text{visc}}}.
\end{equation}

We are interested in the steady states, i.e., solutions where $\frac{\partial \mathbf{u}}{\partial t}$ vanishes.
The high-order simulations are computed through time-advancement, while the RB reduced order solutions are computed 
through fixed-point iterations of a nonlinear solver.

\subsection{Nonlinear solver}

The \emph{Oseen}-iteration is a secant modulus fixed-point iteration, which in general exhibits a linear rate of convergence.
It solves for a steady-state solution, i.e., $\frac{\partial \mathbf{u}}{\partial t} = 0$ is assumed.
Given a current iterate (or initial condition) $\mathbf{u}^k$, the next iterate $\mathbf{u}^{k+1}$ is found
by solving the following linear system:
\begin{eqnarray}
 -\nu_{\text{visc}} \Delta \mathbf{u}^{k+1} + (\mathbf{u}^k \cdot \nabla) \mathbf{u}^{k+1} + \nabla p &=& \mathbf{f}  \text{ in } \Omega, \label{Hess:eq_Oseen_main} \cl
\nabla \cdot \mathbf{u}^{k+1} &=& 0   \text{ in } \Omega, \cl
 %\mathbf{u} &=& f_d  \text{ on } \partial \Omega,
\mathbf{u}^{k+1} &=& \mathbf{d}  \text{ on } \Gamma_D, \cl
\nabla \mathbf{u}^{k+1} \cdot \mathbf{n} &=& \mathbf{g} \text{ on } \Gamma_N. \el
\end{eqnarray}
Iterations are stopped when the relative difference between iterates falls below a predefined tolerance in a suitable norm, like the $L^2(\Omega)$ or $H^1_0(\Omega)$ norm. 

\section{Sparse Polynomial Interpolation}

The presented sparse polynomial interpolation approach is based on the literature references \cite{10.1007/s10208-013-9154-z} and \cite{CHKIFA2015400}.
Let $\mathcal{P}$ denote a parameter domain of vectors $\mathbf{y} = (y_1, \ldots, y_d) \in \mathbb{R}^d$ with $d$ the number of parameters.
Each parameter direction has been normalized to the interval $\left[ -1, 1 \right] $, such that

\begin{equation}
 \mathcal{P} = \left[ -1, 1 \right]^d \subset \mathbb{R}^d .
\end{equation}

Introduce the parameter-to-solution map as 

\begin{equation}
 \mathbf{y} \in \mathcal{P} \mapsto \mathbf{u}(\mathbf{y}) \in X ,
\end{equation}

\noindent assuming the well-posedness of \eqref{Hess:NSE0}-\eqref{Hess:NSE1} for each $\mathbf{y} \in \mathcal{P}$ and $X$ denotes a suitable function space, such as $H^1_0(\Omega)$ with model-specific boundary conditions.

Let $\mathcal{F}$ denote the set of all finitely supported sequences $\nu = (\nu_1, \nu_2, \ldots, 0, 0, \ldots) \in \mathbb{N}^\mathbb{N}_0$.
Since this work only considers models with a finite number of parameters, $\mathcal{F}$ will be set to $\mathbb{N}^d_0$.
Introduce the ansatz of expressing $\mathbf{u}(\mathbf{y})$ as a high-order polynomial approximation

\begin{equation}
 \mathbf{u}(\mathbf{y}) = \sum_{\nu \in \Lambda} \mathbf{u}_\nu \mathbf{y}^\nu ,
\end{equation}

\noindent with $\Lambda \subset \mathcal{F}$ a finite set of multiindices, $\mathbf{u}_\nu$ is a solution to \eqref{Hess:NSE0}-\eqref{Hess:NSE1} and $\mathbf{y}^\nu = \prod_{j \geq 1} y_j^{\nu_j}$
a polynomial of degree $\sum_{j \geq 1} \nu_j$ over the parameter domain $\mathcal{P}$. The monomials $y_j^{\nu_j}$ are typically not a good choice numerically. Instead, Lagrange interpolation will be used.

\subsection{Univariate Interpolation}

Let $(z_k)_{k \geq 0}$ be a sequence of mutually distinct points in $\left[ -1, 1 \right]$. Let $I_k$ denote the univariate polynomial interpolation operator associated with the first $k$ points of $(z_k)_{k \geq 0}$.
Then the interpolation operator acts on a function $g$, which is defined over the parameter domain $\mathcal{P}$. 
The interpolation operator is defined as 

\begin{equation}
 I_k g = \sum_{i = 0}^k g(z_i) l_i^k ,
\end{equation}

\noindent with the Lagrange polynomials

\begin{equation}
 l_i^k = \prod_{j = 0, j \neq i}^k \frac{y - z_j}{z_i - z_j} .
\end{equation}

Introduce the difference operator 

\begin{equation}
 \Delta_k = I_k - I_{k-1} ,
\end{equation}

\noindent where $I_{-1}$ is the null operator. As a consequence $\Delta_0 g = I_0 g = g(z_0)$, the constant polynomial with value $g(z_0)$.

Thus, it holds that 

\begin{equation}
 I_n = \sum_{k=0}^n \Delta_k .
\end{equation}

Introduce the hierarchical polynomials of degree $k$ as

\begin{equation}
 h_k(y) = \prod_{j=0}^{k-1} \frac{y - z_j}{z_k - z_j}, \quad k > 0 , \quad h_0(y) = 1 ,
\end{equation}

\noindent which implies 

\begin{equation}
 \Delta_k g = \alpha_k(g) h_k ,
\end{equation}

\noindent with 

\begin{equation}
 \alpha_k(g) = g(z_k) - I_{k-1} g(z_k) .
\end{equation}

This allows the representation

\begin{equation}
 I_n g = \sum_{k=0}^n \alpha_k(g)h_k .
 \label{univ_interpolation}
\end{equation}

\subsection{Tensorization}

In the multivariate case consider the parameter vector $\mathbf{y} = (y_1, \ldots, y_d) \in \mathbb{R}^d$.
Given a multiindex $\mathbf{\nu} \in \mathcal{F}$, define the multivariate point $z_\mathbf{\nu}$ as

\begin{equation}
 z_\nu = (z_{\nu_j})_{j \geq 1} \in \mathcal{P} .
 \label{multivariate_point}
\end{equation}

The tensorized hierarchical function is 

\begin{equation}
 H_\mathbf{\nu}(\mathbf{y}) = \prod_{j \geq 1} h_{\nu_j}(y_j) ,
\end{equation}

\noindent and the tensorized multivariate operators are 

\begin{equation}
 I_\mathbf{\nu} = \bigotimes_{j \geq 1} I_{\nu_j} , \quad \Delta_\mathbf{\nu} = \bigotimes_{j \geq 1} \Delta_{\nu_j} .
\end{equation}

\subsection{Sparse Interpolation Operator}

For two multiindices $\nu$ and $\mu$, the relation $\nu \geq \mu$ is defined as $\nu_i \geq \mu_i$ for all parameter directions $i$.

An index set $\Lambda$ of (multi-)indices $\nu$, which fulfills the property that 

\begin{equation}
\left( \nu \in \Lambda \text{ and } \nu \geq \mu \right) \Rightarrow \mu \in \Lambda,
\end{equation}

\noindent is called a \emph{monotone} or \emph{downward closed} set.

Given a downward closed set $\Lambda$, the sparse interpolation operator is defined as

\begin{equation}
 I_\Lambda = \sum_{\nu \in \Lambda} \Delta_\nu .
 \label{sparse_int_op}
\end{equation}

The efficient hierarchical computation of the sparse interpolation operator is shown in \cite{10.1007/s10208-013-9154-z}.

\subsection{Leja Points}

The suggested point rules in \cite{10.1007/s10208-013-9154-z} are Leja sequences, composed of Leja points.
Leja points are defined recursively by maximizing 

\begin{equation}
 F^N(y) = \prod_{i=1}^{N-1} \vert (y-x_i) \vert
\end{equation}

\noindent over $\left[ -1, 1 \right] $ for a given initial $x_1$, such that

\begin{equation}
 x_N = \arg \max_{y \in \left[ -1, 1 \right]} F^N(y) .
 \label{leja_max}
\end{equation}

Solving \eqref{leja_max} is actually computationally hard, such that in practice a fine grid of the interval $\left[ -1, 1 \right]$ is used.
Higher dimensional Leja points are then determined by tensorization and in particular \eqref{leja_max} is not re-expressed in higher dimensions.
Symmetrized Leja points are defined by choosing $x_1 = 0, x_2 = 1, x_3 = -1$ and then evaluating \eqref{leja_max} for even $N$ while choosing $x_N = x_{N-1}$ for odd $N$.
A set of points can be put in Leja ordering by restricting the maximization in \eqref{leja_max} to the set itself.

%\section{theoretical considerations}

\section{Numerical Simulations}

The sparse polynomial approach is used to generate reduced order models for parametrized channel flows with one and two parameters.
To access the quality of the approximations, they are compared against reduced basis (RB) methods based on the proper orthogonal decomposition (POD).

The models are discretized with the spectral element method (SEM) \cite{karniadakis1999spectral} using the PDE framework \texttt{Nektar++}\footnote{ \texttt{www.Nektar.info} } and the model reduction software 
\texttt{ITHACA-SEM}\footnote{ \texttt{https://github.com/mathLab/ITHACA-SEM} }.

\subsection{Channel with a Narrowing of Varying Width}

\begin{figure}
 \includegraphics[scale=.2]{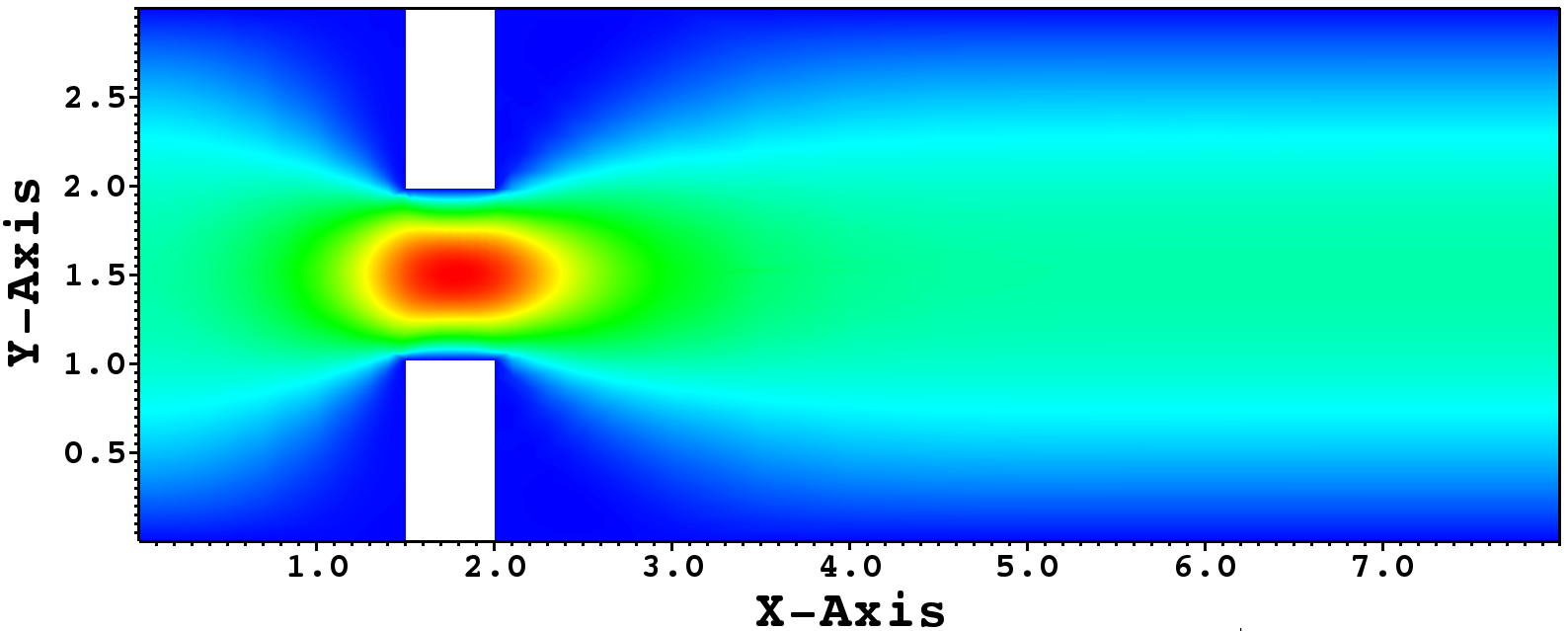} $\quad$
 \includegraphics[scale=.35]{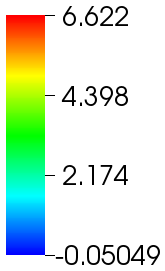}  \\
 \includegraphics[scale=.2]{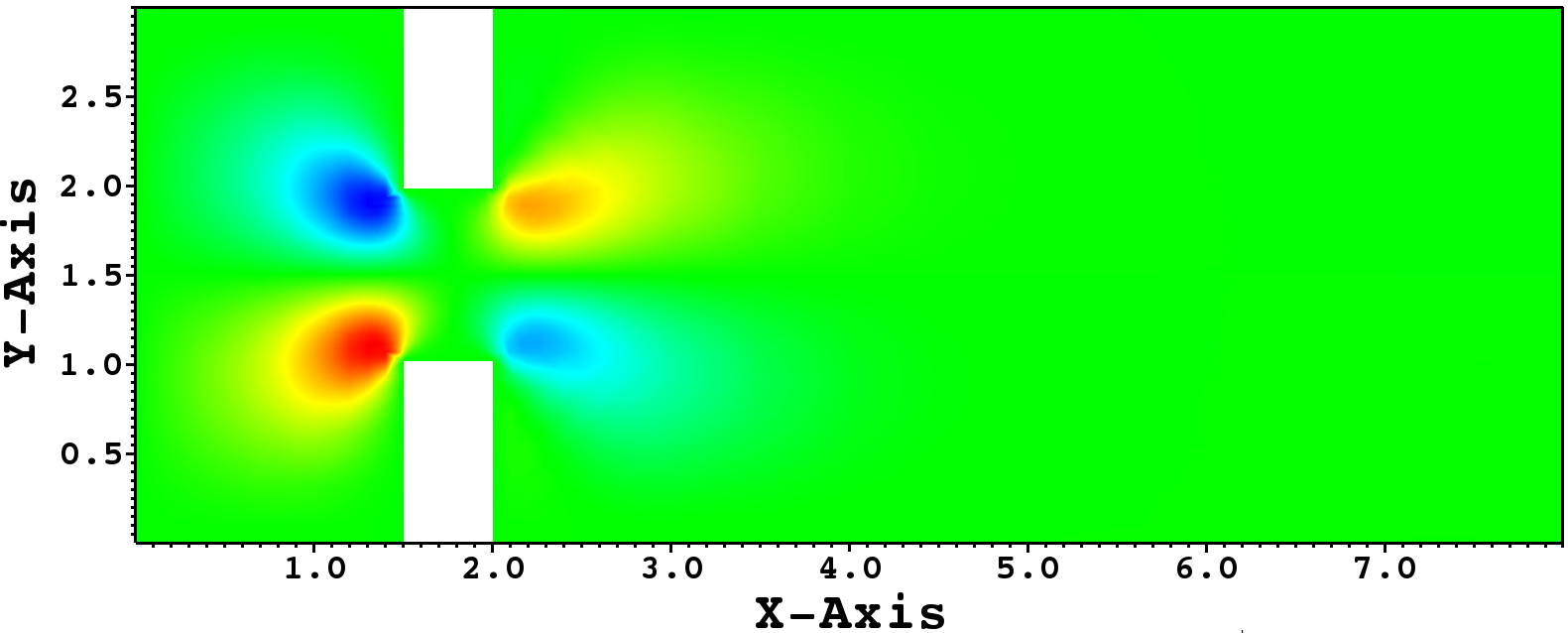} $\quad$
 \includegraphics[scale=.35]{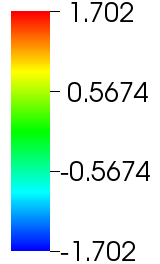}
 \caption{Full order, steady-state solution for $\mu = 1$: velocity in x-direction (top) and y-direction (bottom).}
 \label{Hess:Geo1}
\end{figure}

Consider a channel flow with a narrowing of varying width. The velocity field solution at the reference parameter $\mu_{ref} = 1$ is shown in Fig.~\ref{Hess:Geo1}.
Some more field solutions are shown in \cite{10.1007/978-3-030-39647-3_45} and closely related models have been computed also in \cite{Hess2019CMAME}, \cite{HessRozza2019} and \cite{PichiHesthaven2021}.
The geometry is decomposed into $36$ triangular spectral elements and the velocity is resolved with modal Legendre polynomials of order $11$.
The inflow profile on the left side is parabolic with $u_x(0,y) = y(3-y)$ for $y \in [0, 3]$.
At the outlet, a stress-free boundary condition is set and everywhere else hold no-slip conditions.

The parameter domain is $\mu \in \left[ 0.1 , 2.9 \right]$, which is affinely mapped to the interval $\nu \in \left[ -1 , 1 \right]$ to conform with 
the sparse polynomial approximation assumptions.
With a changing parameter, the geometry always remains symmetric to the horizontal centerline at $y = 1.5$.
The kinematic viscosity is kept constant at $\nu_{\text{visc}} = 1$.

\begin{figure}
\begin{center}
\includegraphics[width=4in]{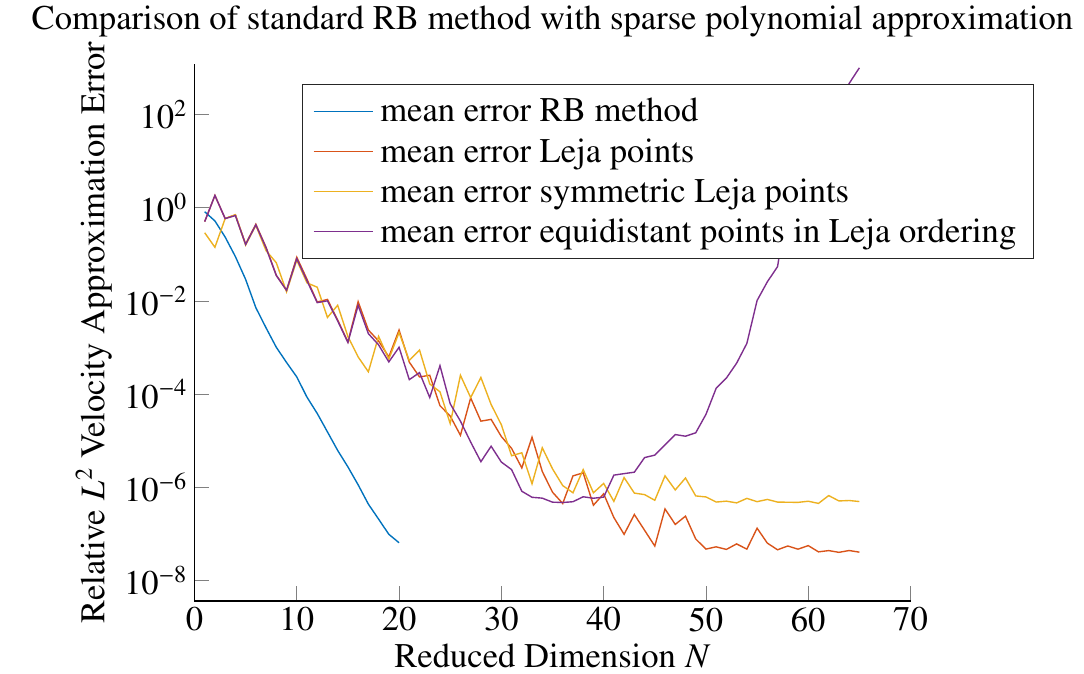} 
\caption{Comparison of mean error of the RB method with sparse polynomial approximations using Lagrange polynomials and various point rules.} 
\label{fig:model1_point_rules_mean_error}
\end{center}
\end{figure}

\begin{figure}
\begin{center}
\includegraphics[width=4in]{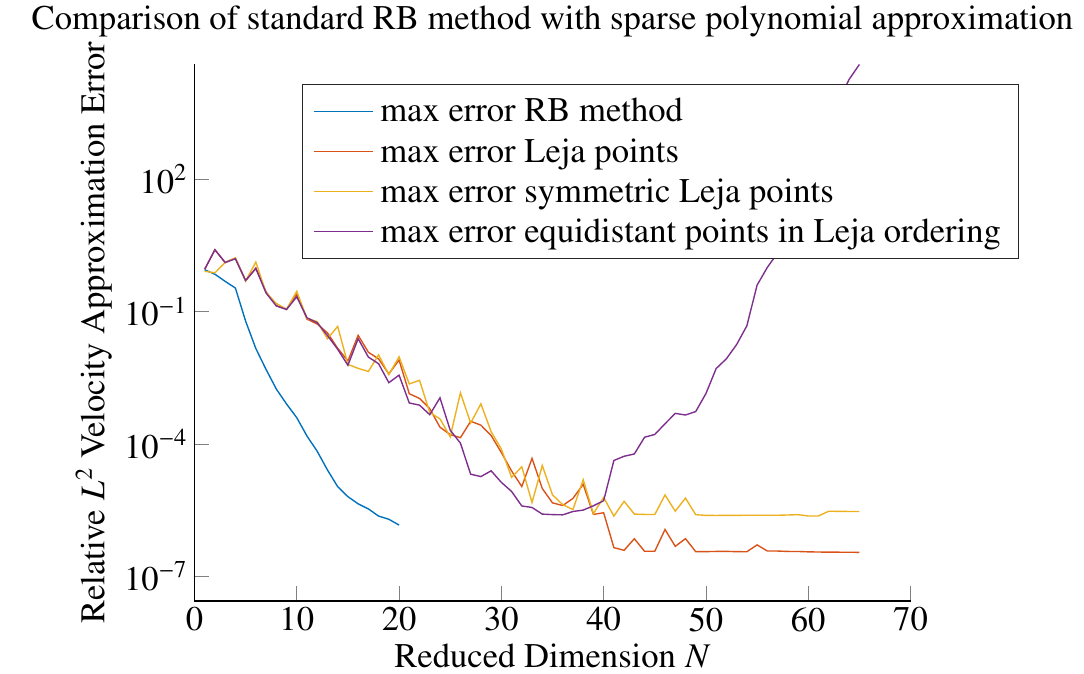} 
\caption{Comparison of maximum error of the RB method with sparse polynomial approximations using Lagrange polynomials and various point rules.} 
\label{fig:model1_point_rules_max_error}
\end{center}
\end{figure}

The parametric variation in geometry allows an affine decomposition of the Navier-Stokes element matrix in the parameter.
In particular, Fig.~\ref{fig:model1_point_rules_mean_error} and Fig.~\ref{fig:model1_point_rules_max_error} show the mean and maximum error 
with growing reduced order dimension for the standard reduced basis procedure, the Lagrange polynomials with Leja points, 
the Lagrange polynomials with symmetrized Leja points and Lagrange polynomials with equidistant points in Leja ordering.

Although the sparse polynomial interpolation and the RB are two very different approximation algorithms, the results are both plotted versus the reduced dimension $N$.
In both cases the size of the ROM model is compared. 
For the sparse interpolation the reduced dimension refers to the number of PDE solutions in the polynomial expansion, i.e., the $n$ in \eqref{univ_interpolation}.
For the RB the reduced dimension refers to the number of POD modes used to project the equations. 
This is not a perfectly fair comparison, since more PDE solutions were used in the RB method to generate the POD modes (namely $40$ for this model), but it still allows to draw some conclusions.

It can be observed that all methods show a faster than polynomial order convergence, since the slope is linear in a semi-log plot. 
The three sparse polynomial approximations reach a mean accuracy of six digits in the velocity at a ROM dimension of about $35$ and five digits accuracy in the maximum at this dimension.
The RB reaches a mean accuracy of six digits in the velocity at a ROM dimension of $20$ and five digits accuracy in the maximum at ROM dimension $20$.

All three point rules provide a similar approximation quality up a reduced dimension of $35$.
The equidistantly distributed points in Leja ordering then diminish in approximation quality, which is a known phenomenon, since high order Lagrange interpolants without special choice of interpolation points are ill-conditioned.
Using equidistantly distributed points \emph{without} Leja ordering provides no approximation at all and has consistently a mean relative error of about $100\%$.

The approximations with Leja points do not provide a better approximation than single precision, which can usually be improved up to double precision by re-arranging how terms are computed. 
The same holds true for RB approximations. However, a stable approximation with six digits of accuracy is usually enough in practical applications.

\subsection{Channel with a Narrowing of Varying Curvature}

A two parameter model is considered with parametric variation in the curvature of the narrowing and variation in the kinematic viscosity.
This model was analysed in \cite{MR4099821} and the results established there with the RB method and the \emph{empirical interpolation method} (EIM) 
will serve to compare the accuracy of ROMs by sparse polynomial interpolation.

\begin{figure}[ht]
\begin{center}
 \includegraphics[scale=.23]{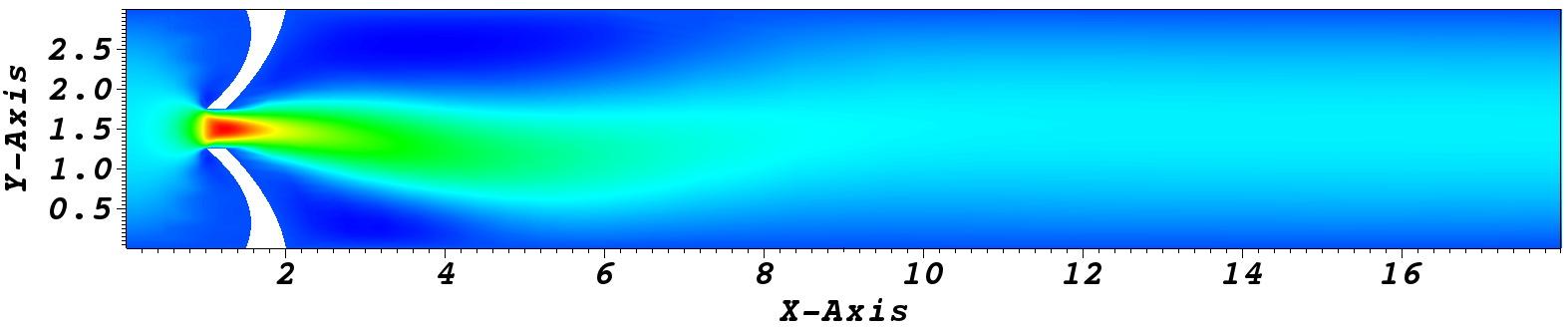} $\quad$
 \includegraphics[scale=.35]{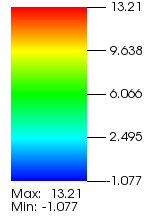} \\
 \includegraphics[scale=.23]{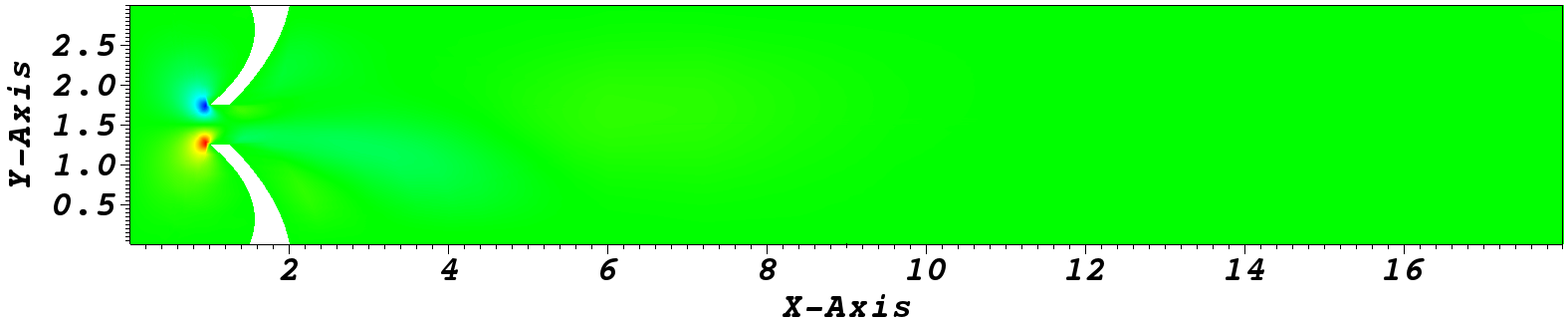} $\quad$
 \includegraphics[scale=.35]{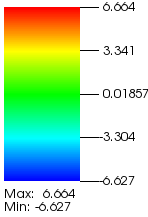}
 \end{center}
 \caption{Full order, steady-state solution in the geometry with curved walls and for $\nu = 0.15$: velocity in x-direction (top) and y-direction (bottom).}
 \label{Hess:FOM_curved_small_nu}
\end{figure}

\begin{figure}[ht]
\begin{center}
 \includegraphics[scale=.23]{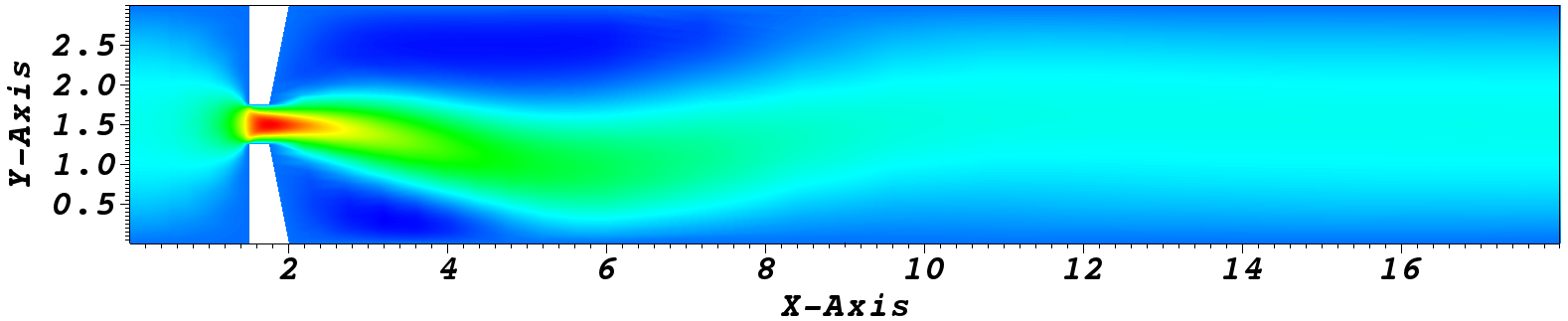} $\quad$
 \includegraphics[scale=.35]{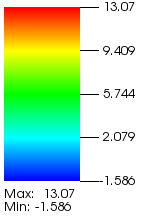} \\
 \includegraphics[scale=.23]{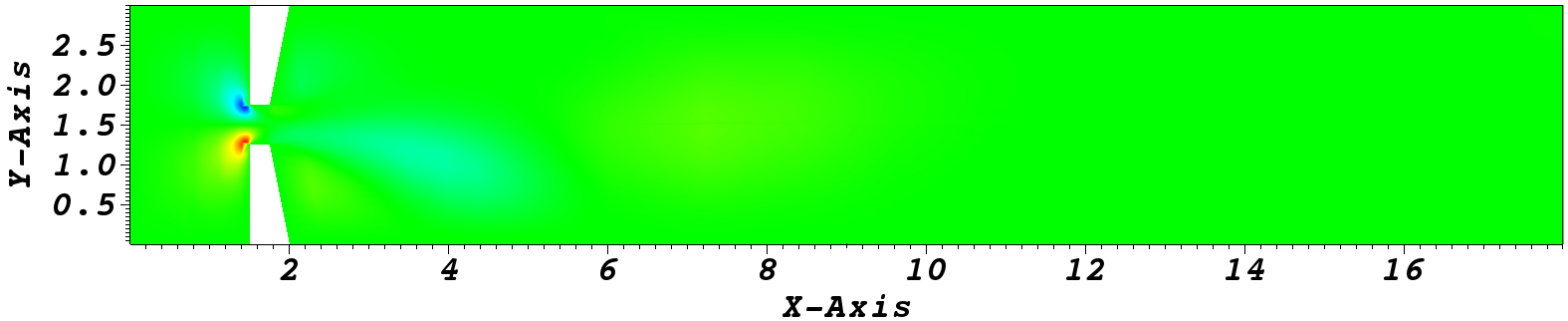} $\quad$
 \includegraphics[scale=.35]{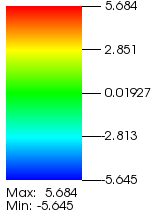}
 \end{center}
 \caption{Full order, steady-state solution in the geometry with straight walls and for $\nu = 0.15$: velocity in x-direction (top) and y-direction (bottom).}
 \label{Hess:FOM_straight_small_nu}
\end{figure}

\begin{figure}[ht]
\begin{center}
 \includegraphics[scale=.23]{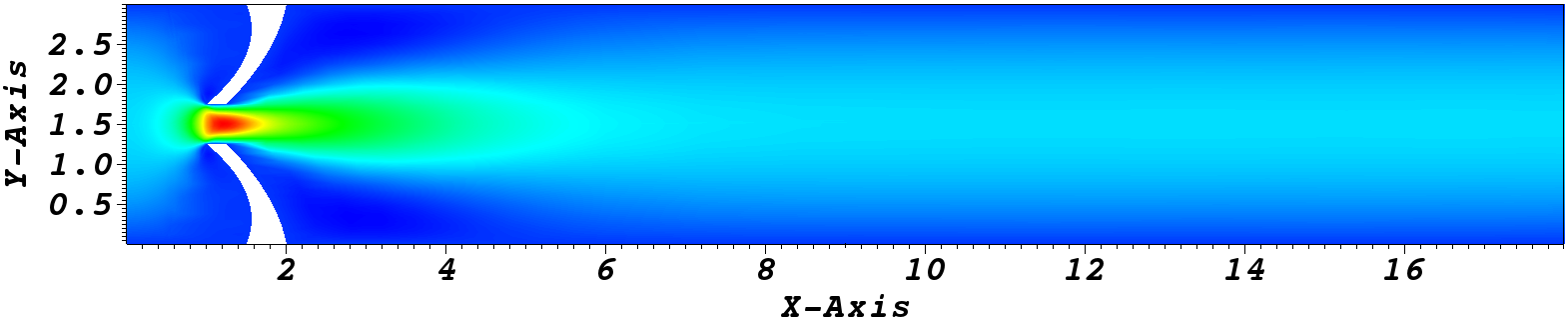} $\quad$
 \includegraphics[scale=.35]{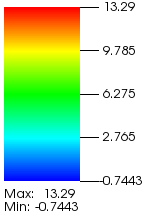} \\
 \includegraphics[scale=.23]{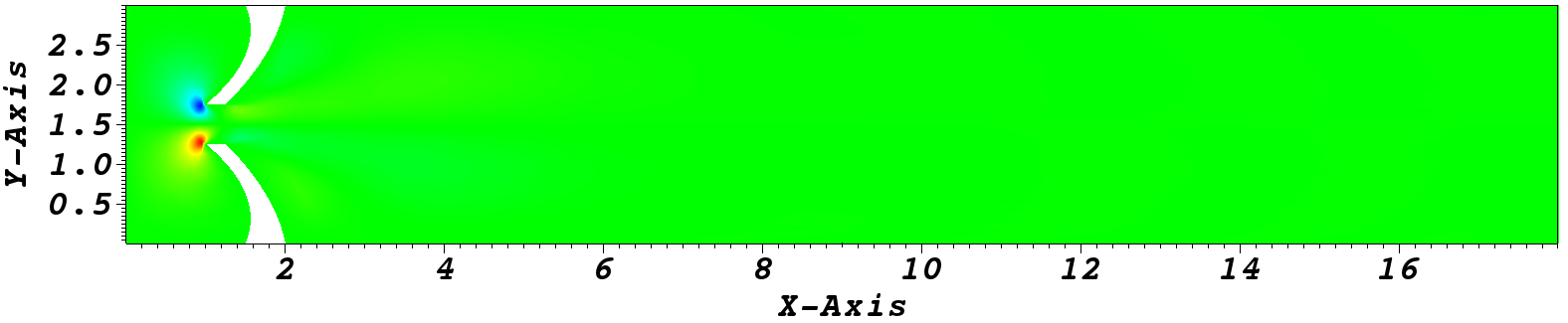} $\quad$
 \includegraphics[scale=.35]{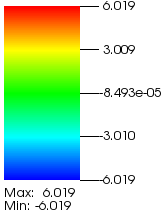}
 \end{center}
 \caption{Full order, steady-state solution in the geometry with curved walls and for $\nu = 0.2$: velocity in x-direction (top) and y-direction (bottom).}
 \label{Hess:FOM_curved_large_nu}
\end{figure}

\begin{figure}[ht]
\begin{center}
 \includegraphics[scale=.23]{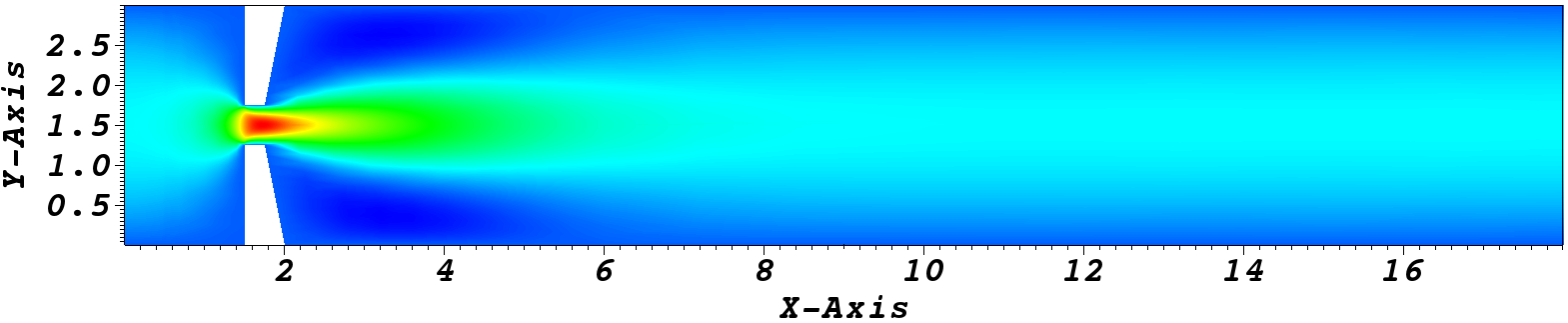} $\quad$
 \includegraphics[scale=.35]{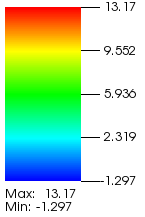} \\
 \includegraphics[scale=.23]{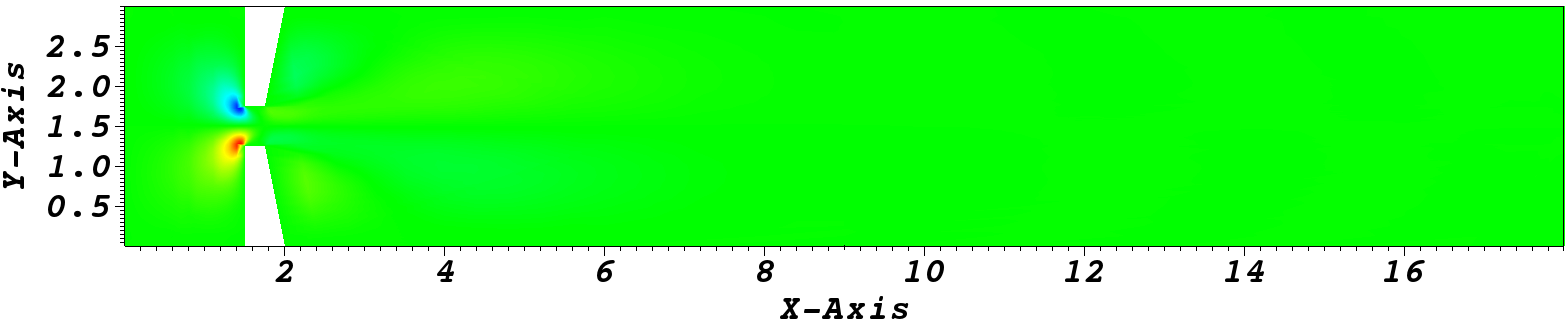} $\quad$
 \includegraphics[scale=.35]{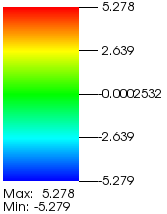}
 \end{center}
 \caption{Full order, steady-state solution in the geometry with straight walls and for $\nu = 0.2$: velocity in x-direction (top) and y-direction (bottom).}
 \label{Hess:FOM_straight_large_nu}
\end{figure}

Consider the channel flow through a narrowing created by walls of varying curvature and with variable kinematic viscosity.
See Fig.~\ref{Hess:FOM_curved_small_nu} and Fig.~\ref{Hess:FOM_straight_small_nu} for the steady-state velocity components
for $\nu = 0.15$ in a geometry with curved walls and straight walls, respectively.
Fig.~\ref{Hess:FOM_curved_large_nu} and Fig.~\ref{Hess:FOM_straight_large_nu} show the steady-state velocity components
for $\nu = 0.2$ in a geometry with curved walls and straight walls, respectively. 
These are the four corners of the rectangular parameter domain and constitute the most extreme solutions. 
Fig.~\ref{Hess:FOM_curved_small_nu} and Fig.~\ref{Hess:FOM_curved_large_nu} also show the strongest curvature of all configurations.

The spectral element expansion uses modal Legendre polynomials of order $p = 10$ for the velocity. 
The pressure \emph{ansatz} space is chosen of order $p-2$ to fulfill the inf-sup stability condition (\cite{B-quarteroniv2, Brezzi:mixed}).
A parabolic inflow profile is  prescribed at the inlet (i.e., $x = 0$) with horizontal velocity component  $u_x(0,y) = y(3-y)$ for $y \in [0, 3]$.
At the outlet (i.e., $x = 18$) a stress-free boundary condition is imposed, while everywhere else a no-slip condition is prescribed.
Symmetric boundary conditions are considered in order to study
the symmetry breaking due to the nonlinearity in problem \eqref{Hess:NSE0}--\eqref{Hess:NSE1}.
A more realistic setting considers also different inlet velocity profiles and the pulsatility of the flow and would then include the Strouhal number as a parameter.

The viscosity varies in the interval $\nu_{\text{visc}} \in [0.15, 0.2]$.
The Reynolds number $Re$, see Eq.~\eqref{eq:re}, depends on the kinematic viscosity. 
As $Re$ is varied for each fixed geometry, a supercritical pitchfork bifurcation occurs:
for $Re$ higher than the critical bifurcation point, three solutions exist.
Two of these solutions are stable, one with a jet towards the top wall and one with a jet
towards the bottom wall and one is unstable. The unstable solution 
is symmetric to the horizontal centerline at $y=1.5$, while the jet of the stable solutions
 undergoes the \emph{Coanda effect}.

In this investigation, we do not deal with recovering all bifurcation branches, but limit our attention 
to the stable branch of solutions with jets hugging the bottom wall. However, we remark that
recovering all bifurcating solutions with model reduction methods is also possible, see, e.g., (\cite{Hess:Herrero2013132}) and (\cite{PintorePichiHessRozzaCanuto2021}).

\subsubsection{Generating Curved Geometries}

The different curvatures are approximated by polynomials. 
Each curved wall is defined by a second order polynomial, interpolating three prescribed points. 
While the points at the domain boundary $y=0$ and $y=3$ are kept fixed, 
the inner points are moved towards $x=0$ in order to create an increasing curvature.
The tip of both narrowings and one intermediate point between the tip and the wall are prescribed, while the point where the wall and the narrowing meet remains constant.
These three points define a quadratic polynomial, which is used to model the edge from the tips to upper and lower wall, respectively.
This is a standard feature implemented in the PDE solver \texttt{Nektar++}.

The mesh remains topologically equivalent for each parametric configuration. This allows to easily map the mesh to the reference configuration via a plain pullback (see \cite{MR3780742} for a discussion of the plain pullback) 
acting directly on the degrees of freedom, i.e., the entries of the solution vector.

\subsubsection{Numerical Results}

\begin{figure}[ht]
\begin{center}
 \includegraphics[scale=.2]{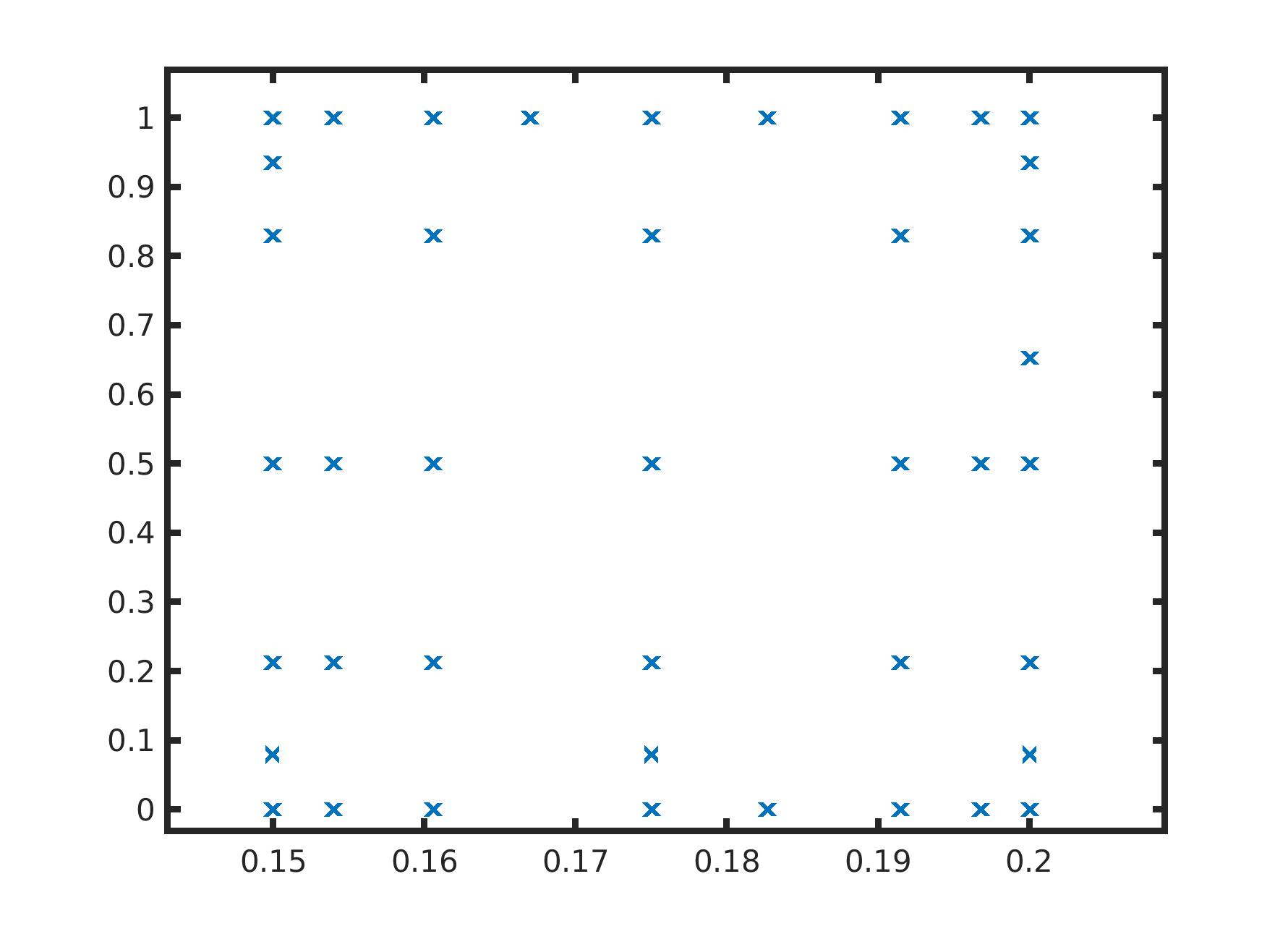} 
\end{center}
 \caption{Chosen Leja points plotted with the kinematic viscosity on the x-axis and a measure of the curvature on the y-axis, where the `1' refers to the maximum curvature.}
 \label{Hess:leja_2d}
\end{figure}

The Leja points are computed as in Eq.~\eqref{leja_max} in each parameter direction and a tensorized grid of the two dimensional parameter domain is used.
A grid rule has to be chosen that determines how the sequence of multiindices $(\nu^n)_{n \geq 1}$ has to be chosen that defines the multivariate points $z_{\nu^n}$ of Eq.~\eqref{multivariate_point}.
Here, a downward closed set is formed in a canonical way by choosing $\nu^1 = (0, 0)$, $\nu^2 = (1, 0)$, $\nu^3 = (0, 1)$, $\nu^4 = (2, 0)$, $\nu^5 = (1, 1)$, \ldots, i.e., increasing 
the sum $\nu^n_1 + \nu^n_2$ only after all possible combinations of elements with the same sum $\nu^n_1 + \nu^n_2$ have been added to $(\nu^n)_{n \geq 1}$.
The chosen points are depicted in Fig.~\ref{Hess:leja_2d}.

The set $\Lambda$, which defines the sparse interpolation operator Eq.~\eqref{sparse_int_op} is given by $\Lambda_n = \{ \nu^i, i = 1, \ldots, n  \}$.
This defines a hierarchical sequence of index sets $\Lambda_1 \subset \Lambda_2 \subset \Lambda_3 \subset, \ldots, \Lambda_n$, which allows to reuse the computed snapshot solutions when
updating the interpolation operator $I_{\Lambda_{n-1}}$ to $I_{\Lambda_n}$.

\begin{figure}[ht]
\begin{center}
 \includegraphics[scale=1]{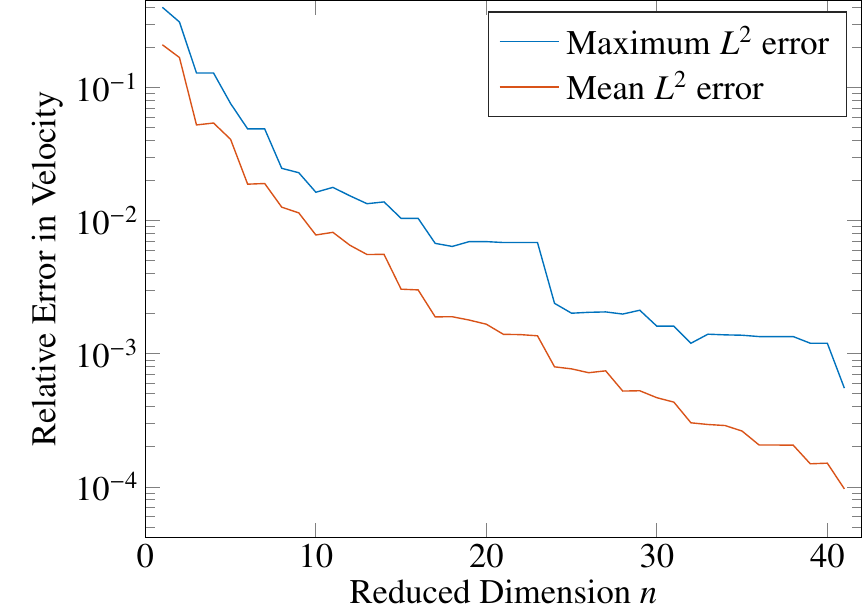} 
\end{center}
 \caption{Relative error for increasing sparse interpolation ROM dimension evaluated over a fine grid of $72$ snapshots.}
 \label{Hess:Num2_rel_L2}
\end{figure}

The interpolation operators are computed up to dimension $41$ and the maximum and mean $L^2$ errors in the velocity are computed for $72$ reference solutions. 
Only the Leja-points without explicit symmetrization are used, since the first numerical test could show that the results of the point rules are similar.
Fig.~\ref{Hess:Num2_rel_L2} shows the relative error in the velocity for increasing ROM size of the sparse interpolation.
A maximum error of less than $1\%$ is reached at dimension $17$ and a maximum error of less than $0.1\%$ is reached at dimension $41$.
A mean error of less than $1\%$ is reached at dimension $10$ and a mean error of less than $0.1\%$ is reached at dimension $24$.
The error does not jump above these thresholds for higher dimensions.
This indicates, that the sparse polynomial interpolation generates usable and reliable ROMs, that can be refined to higher accuracy as needed.

The same model has been investigated using a reduced-basis (RB) framework with the empirical interpolation method (EIM) in \cite{MR4099821}.
The reduced order model with $N = 20$ basis functions, showed an absolute error at the bifurcation point value of less than $0.01$ at $46$ parameter locations and less than $0.1$ at $63$ parameter locations.
The chosen  bifurcation point is the vertical velocity at the point (2, 1.5), which was used in \cite{MR4099821} to access the accuracy.
This means, that the RB method was not able to generate overall accurate methods in each of the $72$ test points, in contrast to the sparse polynomial interpolation.
To plot the RB approximation accuracy as in Fig.~\ref{Hess:Num2_rel_L2} would not add helpful information as the error was jumping up frequently when increasing the ROM dimension.

Another issue is that the EIM relies on the fast computation of a few matrix entries during the online phase. 
However, the spectral element \emph{ansatz} functions have support over a whole spectral element, so this operation cannot be performed as fast as 
with a finite element method for example. The speed-up will thus not be as significant when using EIM.
Also in the case of a finite-element or finite-volume discretization the sparse polynomial approach would likely perform better, since the sparse polynomial approach avoids the EIM completely. 
However, the gain in computational speed-up would not be as significant when comparing to a finite-element or finite-volume discretization.
\newline

A few points are worth highlighting:
\newline

\textit{Relevance of Results}
\newline

The sparse polynomial interpolation was able to compute accurate ROMs, which are stable when increasing the ROM size.
While the first numerical example shows that the ROM sizes are not as small as for the RB method, the RB still needs to compute a similar number of overall snapshots for the POD sample matrix.
Taking this into account, both methods show a comparable performance.
The curved boundaries in the second numerical example allow topologically equivalent mesh changes, but introduce a parametric nonlinearity.
The sparse interpolation nevertheless produced accurate ROMs, where the RB largely failed.
Parametric nonlinearities are inherently difficult for many ROM methods, but at least in the case of topologically equivalent meshes, the sparse polynomial interpolation could be a method of choice.
Additionally, the sparse polynomial interpolation offers adaptivity in choosing the interpolation points with a heuristic error estimator, see \cite{10.1007/s10208-013-9154-z}.
\newline

\textit{Offline-Online Decomposition}
\newline

The offline-online decomposition separates the computations in two parts. 
The offline phase performs time-intensive computations such as the snapshot computations, while the online phase is quickly solves the ROM for many parameters of interest or in a real-time context.
The offline-online splitting is also present in the sparse interpolation, since the snapshot solutions can be computed on a high-performance cluster (HPC), while the evaluation of the sparse
operator for a parameter of interest can be done efficiently without a large computational effort. This is important for a wide applicability of the method and the sparse interpolation shares this property with the RB.
\newline

\textit{Run Time}
\newline

Both methods require the computation of the snapshot solutions. 
There is hardly any additional run time effort for the sparse interpolation, since the interpolation operators can be computed in a hierarchical way as long as the index sets are hierarchical.
The RB on the other hand is much more involved. 
For the assembly of the reduced order systems a reduced trilinear form is computed in the incompressible Navier-Stokes case, which takes approximately as much time as computing the initial snapshot solutions in the shown examples.
Additionally, the compute time for an EIM can be significant in the RB method.
\newline

\textit{Implementation}
\newline

The sparse interpolation can be implemented as outlined in \cite{10.1007/s10208-013-9154-z} with a hierarchical computation.
There are several choices for point rules and polynomials, but the implementational effort is light when compared to the reduced basis method.
The RB for incompressible Navier-Stokes requires to compute the reduced operators of the affine form and the EIM requires to identify the degrees of freedom, which most significantly contribute to the 
system matrix. This is a significant implementation effort, which is not necessary for and has no counterpart in the sparse interpolation.
\newline

\section{Conclusion and Outlook}

The sparse polynomial interpolation generates comparable reduced order models (ROMs) to the reduced basis (RB) method in terms of accuracy and model size.
In terms of applicability, some parametric nonlinearities of the geometry can be treated without altering the method, in particular if the mesh topology remains intact.
Regarding the run time of the method and the ease of implementation, the sparse interpolation is even superior to the RB.
The offline-online splitting is also present in the sparse interpolation, which allows to offload time-consuming snapshot computations to a high-performance cluster, while evaluating the ROM on nearly any machine.
The literature on sparse interpolation offers a lot more than what is discussed here. Namely bounds on the approximation error and techniques for dealing with more complicated nonlinearities as well as adaptive choices of sample sets.
Connecting these topics with the numerical models can be the topic of future research.

\section*{Acknowledgements}

We acknowledge the support provided by the European Research Council Executive Agency by the Consolidator Grant project 
AROMA-CFD ``Advanced Reduced Order Methods with Applications in Computational Fluid Dynamics" - GA 681447, H2020-ERC CoG 2015 AROMA-CFD, PI G. Rozza, and INdAM-GNCS 2019-2020 projects.

\bibliographystyle{plain}

\bibliography{sparsePoly_bib,rbsissa,latexbi}



\end{document}